\documentclass[a4paper, 12pt]{amsart}
\usepackage{amsmath}
\usepackage{amsthm}
\usepackage{amscd}
\usepackage{amssymb}
\usepackage{amsfonts}
\usepackage{latexsym}
\usepackage[mathscr]{eucal}
\usepackage{amsmath,amssymb }

\newtheorem{thm}{Theorem}[section]
\newtheorem{prop}[thm]{Proposition}
\newtheorem{lemma}[thm]{Lemma}

\theoremstyle{remark}
\newtheorem{remark}[thm]{Remark}
\theoremstyle{definition}
\newtheorem{example}[thm]{Example}

\theoremstyle{definition}
\newtheorem{definition}[thm]{Definition}
\theoremstyle{remark}
\newtheorem{notation}[thm]{Notation}

\begin{document}
\title{$C^{\ast }$-algebras arising from substitutions}
\author{MASARU FUJINO}
\address[Masaru Fujino]{Graduate School of Mathematics, Kyushu University, Hakozaki, 
Fukuoka, 812-8581, Japan}
\email{fujino@math.kyushu-u.ac.jp}

\maketitle
\begin{abstract}
In this paper, we introduce a $C^{\ast }$-algebra associated with a 
primitive substitution. We show that when $\sigma$ is proper, 
the $C^{\ast }$-algebra is simple and purely infinite and contains the associated  
Cuntz-Krieger algebra and the crossed product $C^{\ast }$-algebra of 
the corresponding Cantor minimal system. We calculate the $K$-groups.
\end{abstract}
\section{Introduction}
J. Cuntz constructed in his paper~\cite{Cu} a new $C^{\ast }$-algebra 
naturally associated with the $ax+b$-semigroup over $\mathbb{N}$. 
The $C^{\ast }$-algebra $\mathcal{Q}_{\mathbb{N}}$ is generated by a unitary $u$ and 
isometries $\{s_{n}\}_{n\in \mathbb{N^{\ast }}}$ with range projections
$e_{n}=s_{n}s_{n}^{\ast }$ satisfying the relations
$$
s_{n}s_{m}=s_{nm}, ~s_{n}u=u^{n}s_{n}, ~\sum_{k=0}^{n-1}u^{k}e_{n}
u^{-k}=1
$$
for $n, m\in \mathbb{N}^{\ast }$.

Cuntz showed in~\cite{Cu} $\mathcal{Q}_{\mathbb{N}}$ is simple and purely infinite.
Furthermore he showed $\mathcal{Q}_{\mathbb{N}}$ contains Bunce-Deddens algebra 
and is also generated by Bost-Connes algebra adding one unitary which 
corresponds to addition.\\
 In this paper, motivated by the above 
construction of the $C^{\ast }$-algebra, we introduce a $C^{\ast }$-algebra $B_{\sigma }$
associated with a primitive substitution $\sigma$. 
The $C^{\ast }$-algebra $B_{\sigma }$ is generated by 
the Cuntz-Krieger algebra \cite{Cu-K} $\mathcal{O}_{E_{\sigma }^{(2)}}= 
C^{\ast }(\{s_{f}\}_{f\in E_{\sigma }^{(2)}})$,
 where $E_{\sigma }^{(2)}$ is the set of edges with partial order arising from $\sigma$, 
and a unitary $u$ corresponding to addition (Vershik map) 
satisfying the following relations:
$$
u(\sum_{\text{f is maximal in $E_{\sigma }^{(2)}$}}s_{f})
=(\sum_{\text{f is minimal in $E_{\sigma }^{(2)}$}}s_{f})u
$$
and for $f\in E_{\sigma }^{(2)}$ which is not maximal in $E_{\sigma }^{(2)}$,
$$
us_{f}=s_{g},~\text{where $g$ is the successor of $f$.}
$$ 

Let $f_{1},f_{2},\dots ,f_{k-1}$ be maximal elements 
in $E_{\sigma }^{(2)}$ and $f_{k}$ be not a maximal 
element in $E_{\sigma }^{(2)}$. 
Let $g_{k}$ be the successor of $f_{k}$ and $(g_{1},g_{2},\dots ,g_{k})$ be the unique path 
such that $g_{1},\dots ,g_{k-1}$ are minimal. 
By using the relations above, we can calculate
\begin{align*}
us_{f_{1}}s_{f_{2}}\dots s_{f_{k}}
& = u(\sum_{\text{g is maximal in $E_{\sigma }^{(2)}$}}s_{g})s_{f_{2}}\dots s_{f_{k}}\\
& = (\sum_{\text{g is minimal in $E_{\sigma }^{(2)}$}}s_{g})us_{f_{2}}\dots s_{f_{k}}\\
& = \dots = 
(\sum_{\text{g is minimal in $E_{\sigma }^{(2)}$}}s_{g})^{k-1}s_{g_{k}}
& = s_{g_{1}}s_{g_{2}} \dots s_{g_{k}}.
\end{align*}

We show that when $\sigma$ is proper, $B_{\sigma }$ is simple and purely infinite and contains the crossed product $C^{\ast }$-algebra of the corresponding Cantor minimal system, studied by Durand-Host-Skau\cite{Du}, Giordano-Putnam-Skau\cite{Gi} and Herman-Putnam-Skau\cite{He}. 
Furthermore we show that a certain $B_{\sigma }$ has a unique KMS-state and determine the $K$-groups of $B_{\sigma }$.

\section{$C^{\ast }$-algebras associated with substitutions}
In this paper, we study $C^\ast$-algebras associated with primitive substitutions, 
we introduce the  definition of stationary ordered Bratteli diagrams following after~\cite{Du}.

Let $A$ be an alphabet consisting of at least $2$ letters. 
We denote by $A^{+}$ the set of finite words of $A$.
Then a map $\sigma :A\rightarrow A^{+}$ is said to be a substitution on $A$.
By concatenation, for each $n>0$, the $n$-th iteration $\sigma ^{n}:A\longrightarrow A^{+}$
is again a substitution. For any $a\in A$, let $\kappa (a,b)\in \mathbb{N}$ denote the number of occurrences of $b$ in $\sigma (a)$. A substitution $\sigma$ is said to be primitive if there exists $n >0$ 
such that for every $a$, $b\in A$, $b$ occurs in $\sigma ^{n}(a)$.

Let $\sigma$ be a primitive substitution 
on an alphabet $A$ of $N$ elements.  
Then we construct the stationary ordered Bratteli diagram ($V_{\sigma }$,$E_{\sigma }$,$\leq $ )  associated to a substitution $\sigma$ as follows: 
The vertex set $V_{\sigma }$ is given as a disjoint union 
$V_{\sigma }^{(0)}\cup V_{\sigma }^{(1)}\cup V_{\sigma }^{(2)}\cup \dots $ 
where
$V_{\sigma }^{(0)}$ is one point set $\{x_{0}\}$ and each $V_{\sigma }^{(n)}$ $(n\geq 1)$ consists of 
$N$ points. We denote the element of $V_{\sigma }^{(n)}$ by $x(a,n)$ ($a\in A$). 
i.e. $V_{\sigma }^{(n)}=\{x(a,n)| a\in A \}$ $(n\geq 1)$.
The vertex $x(a,n)$ is called the vertex of label $a$ at level $n$.

The edge set $E_{\sigma }$ is given as a disjoint union 
$E_{\sigma }^{(1)}\cup E_{\sigma }^{(2)}\cup E_{\sigma }^{(3)}\cup \dots$  
such that $s(E_{\sigma }^{(n)})\subset V_{\sigma }^{(n-1)}$
and $r(E_{\sigma }^{(n)})\subset V_{\sigma }^{(n)}$ where $r(e)$ denote the range of the edge $e$ and $s(e)$ the source of $e$.
For each $a\in A$, $E_{\sigma } ^{(1)}$ contains exactly one edge from $x_{0}$ to $x(a,1)$ 
for each $a\in A$. 
Each $E_{\sigma }^{(n)}$ consists of edges from $V_{\sigma }^{(n-1)}$ to $V_{\sigma }^{(n)}$.
For an integer $n>1$ and $a$, $b$ in $A$, the number of edges from $x(b,n-1)$ to $x(a,n)$  is equal to $\kappa (a,b)$.

For the Bratteli diagram $(V_{\sigma }$,$E_{\sigma })$ constructed above,
we have a range map $r:E_{\sigma }\longrightarrow V_{\sigma }$ and a source map 
$s:E_{\sigma }\longrightarrow V_{\sigma }$ so that 
$r(E_{\sigma }^{(n)})\subset V_{\sigma }^{(n)}$ and 
$s(E_{\sigma }^{(n)})\subset V_{\sigma }^{(n-1)}$.

For each $n \geq 1$, we define a partial order on $E_{\sigma }^{(n)}$ as follows:

Two edges $e_{1}, e_{2}\in E_{\sigma }$ are compatible if and only if 
$r(e_{1})=r(e_{2})$.

For each $n\geq 2$ and $a\in V_{\sigma }^{(n)}$, we 
put $E_{\sigma }^{(n)}(a)=\{e\in E_{\sigma }^{(n)}\mid $ the label of $r(e)$ is $a.\}$. 
Let $\sigma (a)=a_{1}a_{2},\dots,a_{m}$, then 
we denote the elements of $E_{\sigma }^{(n)}(a)$ by $e(a,k)$ $(1\leq k\leq m)$ 
 so that the label of $s(e(a,k))$ is $a_{k}$ and $r(e(a,k))=a$.

Let $e_{1}$ and $e_{2}$ be in $E_{\sigma }^{(n)} (n\geq 1)$ with $r(e_{1})=r(e_{2})$ 
and let $a$ be the label of $r(e_{1})$, 
then $e_{1}$ and $e_{2}$ are in $E_{\sigma }^{(n)}(a)$. When we write $e_{1}$ and $e_{2}$ 
as $e(a,i)$ and $e(a,j)$ $(1\leq  i,j \leq m)$, we put $e(a,i)<e(a,j)$ if $i<j$.\\

We basically follow the definition and property of ordered Bratteli diagrams 
in~\cite{Du}, \cite{Gi},~\cite{He} and~\cite{Ho}. 
we introduce some necessary notations, definitions
of ordered Bratteli diagrams.
\begin{notation}
Suppose that $(V_{\sigma },E_{\sigma },\leq )$ is the ordered Bratteli diagram
associated to a primitive substitution $\sigma$.\\
(a)Let $x = x(a,n)$ be in $V_{\sigma }^{(n)}$. Then 
$a$ is said to be the label of $x(a,n)$ and $n$ the level of $x(a,n)$. 
Let $L(x)$ denote the label of $x$ and $\ell (x)$ the level of $x$.\\
(b)Let $x = x(a,n) \in V_{\sigma }^{(n)}$. For $k\in \mathbb{Z}$ with $-(n-1)\leq k $,
let $x_{+k}$ denote the vertex $x(a,n+k)\in V_{\sigma }^{(n+k)}$.
Let $f$ be in $E_{\sigma }^{(n)}$ $(n\geq 1)$, then for $k\in \mathbb{Z}$ with 
$-(n-2)\leq k$, we denote by $f_{+k}$ the edge in $E_{\sigma }^{(n+k)}$ such that  
$r(f_{+k}) = r(f)_{+k}$ and $s(f_{+k}) = s(f)_{+k}$.
\end{notation}
 
\begin{definition}
Let $(V_{\sigma },E_{\sigma },\leq )$ be the stationary ordered Bratteli diagram associated with a primitive substitution $\sigma$. 
Let $\mathcal{P}_{\sigma }$ denote the infinite path space associated with $\sigma$. 
i.e. 
$$
\mathcal{P}_{\sigma }=\{(f^{(1)},f^{(2)},\dots)\mid f^{(k)}\in E_{\sigma }^{(k)},~ r(f^{(k)})=s(f^{(k+1)}) ~\text{for any} ~k\geq 1 \}.
$$ 

Let $\widetilde{V}_{\sigma }=V_{\sigma }\backslash (V_{\sigma }^{(0)}\cup V_{\sigma }^{(1)})$. 
For each $x\in \widetilde{V}_{\sigma }$, let $\mathcal{P}_{\sigma }(x)$ denote
the set of finite pathes connecting between the top vertex $x_{0}$ and $x$. i.e.
\begin{align*}
\mathcal{P}_{\sigma }(x) = 
& \{(f^{(1)},f^{(2)},\dots,f^{(m)})|f^{(k)}\in E_{\sigma }^{(k)}
 ~(1\leq k\leq m),\\
& r(f^{(k)}) = s(f^{(k+1)}) ~(1\leq k\leq m-1), r(f^{(m)}) = x \}.
\end{align*}

Let $\mathcal{P}_{\sigma }^{(n)}=\cup _{x\in V_{\sigma }^{(n)}}P_{\sigma }(x)~(n\geq 2)$  
and $\mathcal{P}_{\sigma }^{\text{fin}}$ denote the set $\cup _{n\geq  2}\mathcal{P}_{\sigma }^{(n)}$. 
$\mathcal{P}_{\sigma }^{\text{fin}}$ is called the finite path space arising from $\sigma$. 
\end{definition}

\begin{notation}
Let $\sigma$ be a primitive substitution. \\
(a)We let $E_{\sigma }^{\text{max}}$ and $E_{\sigma }^{\text{min}}$ denote the maximal and minimal edges in $E_{\sigma }$. i.e. 
$$
E_{\sigma }^{\text{max}}=\{f\in E_{\sigma }\mid f ~\text{is maximal in}~ r^{-1}(r(f)).\}. 
$$

Let $\mathcal{P}_{\sigma }^{\text{max}}$ and $\mathcal{P}_{\sigma }^{\text{min}}$ 
denote the set of maximal and minimal 
infinite pathes in $\mathcal{P}_{\sigma }$.
i.e. 
$$
\mathcal{P}_{\sigma }^{\text{max}} = \{(f^{(1)},f^{(2)},\dots )\in \mathcal{P}_{\sigma }\mid f^{(k)}\in E_{\sigma }^{(k)}\cap E_{\sigma }^{\text{max}} , r(f^{(k)}) = s(f^{(k+1)})~ (1\leq k)\}. 
$$

Let $\mathcal{P}_{\sigma }^{\text{fin,max}}$ and $\mathcal{P}_{\sigma }^{\text{fin,min}}$ denote the set of maximal and  minimal finite pathes. i.e. 
$$
\mathcal{P}_{\sigma }^{\text{fin,max}}=
\{(f^{(1)},f^{(2)},\dots ,f^{(m)})\in \mathcal{P}_{\sigma }^{\text{fin}}\mid f^{(k)}\in E_{\sigma }^{(k)}\cap E_{\sigma }^{\text{max}}~(1\leq k\leq m)\}.
$$
(b)For a finite path $f=(f^{(1)}$,$f^{(2)}$,$\dots$,$f^{(n)})$
$\in \mathcal{P}_{\sigma }^{\text{fin}}$, let $s(f)$ and $r(f)$ denote the vertex 
$r(f^{(1)})\in V_{\sigma }^{(2)}$ and $r(f^{(n)})$ respectively and let $\ell (f)$ denote  the length of finite path $f$. i.e.~$\ell (f)=n$.

For an infinite path $l \in \mathcal{P}_{\sigma }$,  
let $s(l)$, the source of $l$, denote the vertex $r(l^{(1)})\in V_{\sigma }^{(1)}$. \\
(c)For any $f=(f^{(1)},f^{(2)},\dots ,f^{(m)})\in \mathcal{P}_{\sigma }^{(m)}$ and  
$g=(g^{(1)}$,$g^{(2)}$,$\dots$,$g^{(n)})\in \mathcal{P}_{\sigma }^{(n)}$ with  $L(r(f))=L(s(g))$, we define the finite path 
$fg\in \mathcal{P}_{\sigma }^{(m+n-1)}$ by $(f^{(1)}$,$f^{(2)}$,$\dots$,$f^{(m)}$,$g_{+(m-1)}^{(2)}$,
$g_{+(m-1)}^{(3)}$,$\dots$,$g_{+(m-1)}^{(n)})$.\\
(d)For any $f=(f^{(1)}$,$f^{(2)}$,$\dots$,$f^{(m)})\in \mathcal{P}_{\sigma }^{(m)}$ and
$l=(l^{(1)}$,$l^{(2)}$,$\dots )
\in \mathcal{P}_{\sigma }$ with $L(r(f))=L(s(l))$, 
we define the infinite path 
$fl\in \mathcal{P}_{\sigma }$ by $(f^{(1)}$,$f^{(2)}$,$\dots$,
$f^{(m)}$,$l_{+(m-1)}^{(2)}$, 
$l_{+(m-1)}^{(3)}$,$\dots)$.
\end{notation}

\begin{definition}
Let $\sigma$ be a primitive substitution.  
Then $(V_{\sigma },E_{\sigma },\leq )$ is said to be properly ordered if 
$\ ^\#\mathcal{P}_{\sigma }^{\text{min}}=\ ^\#\mathcal{P}_{\sigma }^{\text{max}}=1$. 
Then we write $\mathcal{P}_{\sigma }^{\text{min}}=\{p_{min}\}$ and 
$\mathcal{P}_{\sigma }^{\text{max}}=\{p_{max}\}$. 
Now we topologize the infinite path space $\mathcal{P}_{\sigma }$ by the family of cylinder sets
$\{U(f)\}_{f\in \mathcal{P}_{\sigma }^{\text{fin}}}$ where 
$$
U(f)=\{g\in \mathcal{P}_{\sigma }\mid g^{(k)}=f^{(k)}\in E_{\sigma }^{(k)},~
1\leq k\leq \ell (f)\}.
$$ 

Then each $U(f)$ is also closed and we note that 
$\mathcal{P}_{\sigma }$ is a compact Hausdorff space with a countable basis of 
clopen sets. 
Since $\sigma$ is primitive, the cardinality of $\mathcal{P}_{\sigma }$ is infinity  
and so we note that $\mathcal{P}_{\sigma }$ is a Cantor set.
\end{definition}

\begin{definition}
A primitive substitution $\sigma$ on an alphabet $A$ is proper if there exist 
$ b,c\in A$ such that\\
(1)For each $a\in A$, $b$ is the first letter of $\sigma (a)$.\\
(2)For each $a\in A$, $c$ is the last letter of $\sigma (a)$.
\label{proper}
\end{definition}
It is easily seen that if a substitution $\sigma$ is proper, then  
the stationary ordered Bratteli diagram $(V_{\sigma }, E_{\sigma }, \leq )$ is properly ordered.

We recall the definition of Vershik map in \cite{Ho}.
\begin{definition} 
Let $(V_{\sigma },E_{\sigma },\geq )$ be an ordered Bratteli diagram associated with a  primitive substitution $\sigma$. 
For each $f=(f^{(1)}$,$f^{(2)}$,$\dots$ $)\in \mathcal{P}_{\sigma } \backslash \mathcal{P}_{\sigma }^{\text{max}}$,
there exists $k\geq 1$ such that $f^{(1)}$,$\dots$,$f^{(k-1)}\in E_{\sigma }^{\text{max}}$
and $f^{(k)}\not\in E_{\sigma }^{\text{max}}$. 
Let $g^{(k)}$ be the successor of $f^{(k)}$ in $r^{-1}(r(f^{(k)}))$, then we write
the minimal finite path from $x_{0}$ to $s(g^{(k)})$ as 
$(g^{(1)}$,$g^{(2)}$,$\dots$,$g^{(k-1)})$. 
Then the infinite path $(g^{(1)}$,$g^{(2)}$,$\dots$,
$g^{(k)}$,$f^{(k+1)}$,$f^{(k+2)}$,$\dots)$ 
is said to be the successor of $f$.
We let $\lambda _{\sigma }$ be a map defined by
$$
\lambda _{\sigma }:\mathcal{P}_{\sigma } \backslash \mathcal{P}_{\sigma }^{\text{max}} \rightarrow \mathcal{P}_{\sigma }\backslash \mathcal{P}_{\sigma }^{\text{min}},~~
\lambda _{\sigma }(f)=
\text{the successor of }f. 
$$

The map $\lambda _{\sigma }$ is said to be the Vershik map arising from the ordered Bratteli diagram $(V_{\sigma },E_{\sigma },\leq  )$.

When substitution $\sigma$ is proper, 
we extend the Vershik map $\lambda _{\sigma }$ to the whole infinite path space 
$\mathcal{P}_{\sigma }$ by defining $\lambda _{\sigma }(p_{\text{max}})=p_{\text{min}}$. 
Then $\lambda _{\sigma }:\mathcal{P}_{\sigma }\longrightarrow \mathcal{P}_{\sigma }$ is 
a homeomorphism, thus $(\mathcal{P}_{\sigma }$,$\lambda _{\sigma })$ is a topological 
dynamical system. This dynamical system is called the Bratteli-Vershik system 
associated to the ordered Bratteli diagram. It is known in \cite{Ho} 
that when $\sigma$ is proper, 
$\lambda _{\sigma }$ is minimal. 
\end{definition}

\begin{definition} 
Let $\sigma$ be a primitive substitution. 
For any $x\in \widetilde{V}_{\sigma }$, 
let $p_{\text{min}}(x)$ $($resp. $p_{\text{max}}(x)$ $)$ denote the unique minimal $($resp. maximal$)$ path from the top vertex $x_{0}$ to $x$. 
Let $x\in \widetilde{V}_{\sigma }$, then there is a map:

$$
\mathcal{P}_{\sigma }(x)\backslash \{p_{\text{max}}(x)\}\ni f\longrightarrow 
(\text{the successor of}~f)\in \mathcal{P}_{\sigma }(x)\backslash \{p_{\text{min}}(x)\}.
$$

We let the map be expressed by the symbol $\widehat{\lambda} _{\sigma }$. 
\end{definition}

Let $(V_{\sigma },E_{\sigma },\geq )$ be stationary ordered Bratteli diagram associated to a proper primitive substitution $\sigma$. Then we let $T$ be the unitary operator defined by
$$
T:\ell ^{2}(\mathcal{P}_{\sigma }) \longrightarrow \ell ^{2}(\mathcal{P}_{\sigma }),~~~ 
T(\xi _{\ell })=\xi _{\lambda _{\sigma }(\ell )}.
$$

For each $f\in E_{\sigma }^{(2)}$, we let $V_{f}$ be the operator defined by
$$
V_{f}: \ell ^{2}(\mathcal{P}_{\sigma }) \longrightarrow \ell ^{2}
(\mathcal{P}_{\sigma }),~~ 
V_{f}\xi _{l}=\begin{cases}\xi _{fl} & L(r(f))=L(s(l)))\\
0 & L(r(f))\not=L(s(l)))\end{cases}.
$$

Then $\{V_{f}\}_{f\in E_{\sigma }^{(2)}}$ are partial isometries 
 and the next lemma is known.
\begin{lemma} 
Operators 
$\{V_{f}\}_{f\in E_{\sigma }^{(2)}}$ satisfy the Cuntz-Krieger relations. i.e.
$$
\sum_{f\in E_{\sigma }^{(2)}} V_{f}V_{f}^{\ast }=1, ~~~~
$$
$$
V_{g}^{\ast }V_{g}=\sum_{f\in E_{\sigma }^{(2)},L(s(f))=L(r(g))}V_{f}V_{f}^{\ast }
 ~~~(g\in E_{\sigma }^{(2)}).
$$
\label{Cuntz}
\end{lemma}

\begin{lemma}
(a)For any $f\in E_{\sigma }^{(2)}\backslash E_{\sigma }^{\text{max}}$, 
$V_{\widehat{\lambda }_{\sigma }(f)}=TV_{f}$.\\
(b)
$T(\sum_{f\in E_{\sigma }^{(2)}\cap E_{\sigma }^{\text{max}}}V_{f})=
(\sum_{f\in E_{\sigma }^{(2)}\cap E_{\sigma }^{\text{min}}}V_{f})T$.
\label{substitution}
\end{lemma}
\begin{proof}
(a) For any $l\in \mathcal{P}_{\sigma }$,
$$
V_{\widehat{\lambda }_{\sigma }(f)}\xi _{l}=
\begin{cases}
\xi _{\lambda _{\sigma }(fl)} & \text{if} ~L(r(f))=L(s(l)) \\ 0 & \text{otherwise}
\end{cases}
=TV_{f}\xi _{l}
$$
(b)For any $l\in \mathcal{P}_{\sigma }$,
\begin{align*}
T(\sum_{f\in E_{\sigma }^{(2)}\cap E_{\sigma }^{\text{max}}}V_{f})\xi _{l} 
& = 
\xi _{\lambda _{\sigma }(gl)} ~~(g\in 
\mathcal{P}_{\sigma }^{(2)}\cap \mathcal{P}_{\sigma }^{\text{fin,max}} ~\text{with} 
~L(r(g))=L(s(l)))\\
& = \xi _{h\lambda _{\sigma }(l)} ~~(h\in 
\mathcal{P}_{\sigma }^{(2)}\cap \mathcal{P}_{\sigma }^{\text{fin,min}}
 ~\text{with}~ 
L(r(h))=L(s(\lambda _{\sigma }(l))))\\
& = (\sum_{f\in E_{\sigma }^{(2)}\cap E_{\sigma }^{\text{min}}}V_{f})T\xi _{l}.
\end{align*}
\end{proof}

According to the discussion above, we obtain the following definition.
\begin{definition}
Let $\sigma$ be a proper primitive substitution. Then we define $B_{\sigma }$ to be the
universal unital $C^{\ast }$-algebra generated by a unitary $u$ and partial 
isometries $\{s_{f}\}_{f\in E_{\sigma }^{(2)}}$ such that\\
(a)(Cuntz-Krieger relations)
$$
\sum_{f\in E_{\sigma }^{(2)}}s_{f}s_{f}^{\ast }=1,
$$
$$
s_{g}^{\ast }s_{g}=\sum_{f\in E_{\sigma }^{(2)},L(s(f))=L(r(g))}s_{f}s_{f}^{\ast }
~~(f\in E_{\sigma }^{(2)}).
$$
(b)For any $f\in E_{\sigma }^{(2)}\backslash E_{\sigma }^{\text{max}}$, 
$$
s_{\widehat{\lambda }_{\sigma }(f)}=us_{f}.
$$
(c)
$$
u(\sum_{f\in E_{\sigma }^{(2)}\cap E_{\sigma }^{\text{max}}}s_{f})
=(\sum_{f\in E_{\sigma }^{(2)}\cap E_{\sigma }^{\text{min}}}s_{f})u.
$$
\label{def:definition}
\end{definition}

\begin{example}
In \cite{Cu}, Cuntz defined a simple $C^{\ast }$-algebra $B_{2}$ generated by a unitary $u$ and an isometry $s_{2}$ satisfying the following relations.
$$
~s_{2}u=u^{2}s_{2}, ~~~\sum_{k=0}^{2^{n}-1}u^{k}s_{2}^{n}s_{2}^{\ast n}
u^{-k}=1
$$
for any $n\geq 2$.

The $C^\ast$-algebra $B_{2}$ is isomorphic to $B_{\sigma }$ for a certain substitution $\sigma$.  
Let $\sigma$ be a substitution on an alphabet $A=\{a,b\}$ defined by
$\sigma (a)=ab$ and $\sigma (b)=ab$. 
Then $E_{\sigma }^{(2)}=\{f_{1},\widehat{\lambda }_{\sigma }(f_{1}),f_{2},
\widehat{\lambda }_{\sigma }(f_{2})\}$ where 
$f_{1},f_{2}\in E_{\sigma }^{(2)}\cap E_{\sigma }^{\text{min}}$ 
with $L(r(f_{1}))=$a and  $L(r(f_{2}))=$b. 
We put $s_{f_{1}}=s_{2}^{2}s_{2}^{\ast }$, $s_{\widehat{\lambda }_{\sigma }(f_{1})}=us_{2}^{2}s_{2}^{\ast }$, $s_{f_{2}}=u^{2}s_{2}^{2}s_{2}^{\ast }u^{-1}$ 
and $s_{\widehat{\lambda }_{\sigma }(f_{2})}=u^{3}s_{2}^{2}s_{2}^{\ast }u^{-1}$. 
Then $u$, $\{s_{f}\}_{f\in E_{\sigma }^{(2)}}$ satisfy the relations (a)-(c) in Definition 
\ref{def:definition}. 
According to the Theorem \ref{th:simple} in \S3,  
$B_{\sigma }$ is simple. So $B_{\sigma }$ is isomorphic to $B_{2}$.
\end{example}
\begin{remark}
The $C^{\ast }$-algebra associated with substituion can be defined in general primitive substitution, but we restricted the case when $\sigma$ is proper in this paper.

Let $\sigma$ be a primitive substitution. Then we let $T$ be the operator defined by
$$
T:\ell ^{2}(\mathcal{P}_{\sigma }) \longrightarrow \ell ^{2}(\mathcal{P}_{\sigma }),~~~ 
T\xi _{l}=\begin{cases}\xi _{\lambda _{\sigma }(l)} & l\not\in \mathcal{P}_{\sigma }^{\text{max}} \\ 0 & l\in \mathcal{P}_{\sigma }^{\text{max}}
\end{cases}.
$$

Then the next equation clearly holds.
$$
T^\ast T+\chi _{\mathcal{P}_{\sigma }^{\text{max}}}=TT^{\ast }+\chi _{\mathcal{P}_{\sigma }^{\text{min}}}.
$$

Then $T$ and $\{V_{f}\}_{f\in E_{\sigma }^{(2)}}$ satisfy the 
Lemma \ref{Cuntz} and Lemma \ref{substitution}.
It is known that $C^{\ast }(\{V_{f}\mid f\in E_{\sigma }^{(2)}\})$ 
contains the set of compact operators $K(\ell ^{2}(\mathcal{P}_{\sigma }))$. 
Since $\sigma $ is a primitive substitution, the cardinality of $\mathcal{P}_{\sigma }$ is infinity. 
Let $\pi$ be the natural map 
$C^{\ast }(T, \{V_{f} \}_{f\in E_{\sigma }^{(2)}}) \longrightarrow C^{\ast }(T, \{V_{f}\}_{f\in E_{\sigma }^{(2)}})$
/$K(\ell ^{2}(\mathcal{P}_{\sigma }))$. Then $\pi (T)$ is a unitary and $\{\pi (V_{f})\}_{f\in E_{\sigma }^{(2)}}$ 
are partial isometries satisfying the equations of the 
Lemma \ref{Cuntz} and  
Lemma \ref{substitution}.
When $\sigma$ is primitive,  
we can realize the $C^{\ast }$-algebra $B_{\sigma }$ generated by the quatient above.
\end{remark}
\section{simple and purely infiniteness for $B_{\sigma }$}
\begin{notation}
Let $\sigma$ be a proper primitive substitution on an alphabet $A$.\\
(a)For any $f=(f^{(1)},f^{(2)},\dots ,f^{(n)})\in \mathcal{P}_{\sigma }^{\text{fin}}$, we denote by $s_{f}$ the partial isometry 
$s_{f^{(2)}}s_{f^{(3)}_{-1}}\dots s_{f^{(n)}_{-(n-1)}}$.\\
(b) For each $x\in \widetilde{V}_{\sigma }$, 
let $\nu (x)\in \mathbb{N}$ denote the cardinality of the set $\mathcal{P}_{\sigma }(x)$.
\label{notation}
\end{notation}

\begin{lemma}
Let $\sigma$ be a proper primitive substitution.\\
(a) $s_{f}s_{g}^{\ast }=0$ for any $f, g\in \mathcal{P}_{\sigma }^{\text{fin}}$ with 
$L(r(f))\not=L(r(g))$.\\
(b) $\sum_{f\in E_{\sigma }^{(2)}\cap E_{\sigma }^{\text{max}}}s_{f}s_{f}^{\ast }
=(\sum_{g\in E_{\sigma }^{(2)}\cap E_{\sigma }^{\text{max}}}s_{g})
(\sum_{h\in E_{\sigma }^{(2)}\cap E_{\sigma }^{\text{max}}}s_{h})^{\ast }.$\\
(c)For any $f\in \mathcal{P}_{\sigma }^{\text{fin,min}}$ and $0\leq k\leq \nu (r(f))-1$, 
$$
u^{k}s_{f}=s_{\widehat{\lambda }_{\sigma }^{k}(f)}.
$$
(d)Let $f,g\in \mathcal{P}_{\sigma }^{\text{fin,min}}$. If 
$-(\nu (r(f))-1)\leq k\leq \nu (r(g))-1$, then 
$$
s_{f}^{\ast }u^{k}s_{g}\in \{s_{p}^{\ast }s_{q}\mid p,q\in \mathcal{P}_{\sigma }^{\text{fin}}~\text{with}~r(p)=r(f)~\text{and}~r(q)=r(g).\}.
$$
(e)We write $\mathcal{P}_{\sigma }^{\text{min}}=\{p_{\text{min}}\}$. 
Then we let $\delta _{\sigma }$ be a map defined by
$$
\delta _{\sigma }:\mathcal{P}_{\sigma }^{\text{fin,min}}\longrightarrow \{0,1\}, ~
\delta _{\sigma }(f)=\begin{cases}1 &  \text{if}~ L(r(f))=L(s(p_{\text{min}})) \\ 
0 & \text{otherwise}
\end{cases}.
$$
For any $f\in \mathcal{P}_{\sigma }^{\text{fin,min}}$, 
$$
s_{f}(\sum_{g\in E_{\sigma }^{(2)}\cap E_{\sigma }^{\text{min}}}s_{g})
=\delta _{\sigma }(f)\sum_{h\in \mathcal{P}_{\sigma }^{(\ell (f)+1)}\cap 
\mathcal{P}_{\sigma }^{\text{fin,min}}}s_{h}.
$$
(f)For any $f ,g\in \mathcal{P}_{\sigma }^{\text{fin,min}}$ with $L(r(f))=L(r(g))$, we have 
\begin{align*}
s_{f}s_{g}^{\ast } = 
& \sum_{h\in \mathcal{P}_{\sigma }^{(2)}\backslash \mathcal{P}_{\sigma }^{\text{fin,max}}}s_{f}s_{\widehat{\lambda}_{\sigma } (h)}s_{\widehat{\lambda}_{\sigma } (h)}^{\ast }s_{g}^{\ast }\\
& + \delta _{\sigma }(f)u(\sum_{p\in \mathcal{P}_{\sigma }^{(\ell (f)+1)}
\cap \mathcal{P}_{\sigma }^{\text{fin,max}}}s_{p})
(\sum_{q\in \mathcal{P}_{\sigma }^{(\ell (g)+1)}\cap 
\mathcal{P}_{\sigma }^{\text{{fin,max}}}}s_{q})^{\ast }u^{-1}.
\end{align*}
\label{1}
\end{lemma}
\begin{proof}
According to the Definition \ref{def:definition}, (a)-(e) clearly hold.\\
(f) Let $f,g\in \mathcal{P}_{\sigma }^{\text{fin,min}}$, then we have
$$
s_{f}s_{g}^{\ast }=s_{f}u(\sum_{h\in \mathcal{P}_{\sigma }^{(2)}}s_{h}s_{h}^{\ast })u^{-1}s_{g}^{\ast }
$$
$$
= \sum_{h\in \mathcal{P}_{\sigma }^{(2)}\backslash \mathcal{P}_{\sigma }^{\text{fin,max}}}s_{f}s_{\widehat{\lambda }_{\sigma }(h)}s_{\widehat{\lambda}_{\sigma } (h)}^{\ast }s_{g}^{\ast }
+ s_{f}u(\sum_{h\in \mathcal{P}_{\sigma }^{(2)}\cap \mathcal{P}_{\sigma }^{\text{fin,max}}}s_{h}s_{h}^{\ast })u^{-1}s_{g}^{\ast }.
$$

According to the relation (c) in Definition \ref{def:definition} and Lemma 3.2.(b),(e),  we have
\begin{align*}
& s_{f}u(\sum_{h\in \mathcal{P}_{\sigma }^{(2)}\cap \mathcal{P}_{\sigma }^{\text{fin,max}}}s_{h}
s_{h}^{\ast })u^{-1}s_{g}^{\ast }\\
& = s_{f}u(\sum_{h\in \mathcal{P}_{\sigma }^{(2)}\cap \mathcal{P}_{\sigma }^{\text{fin,max}}}s_{h})
(\sum_{i\in \mathcal{P}_{\sigma }^{(2)}\cap \mathcal{P}_{\sigma }^{\text{fin,max}}}s_{i})^\ast u^{-1}s_{g}^{\ast }\\
& = s_{f}(\sum_{h\in \mathcal{P}_{\sigma }^{(2)}\cap 
\mathcal{P}_{\sigma }^{\text{fin,min}}}s_{h})uu^{-1}
(\sum_{i\in \mathcal{P}_{\sigma }^{(2)}\cap \mathcal{P}_{\sigma }^{\text{fin,min}}}s_{i})^\ast s_{g}^{\ast }\\
& = \delta _{\sigma }(f)(\sum_{h\in \mathcal{P}_{\sigma }^{(2)}\cap 
\mathcal{P}_{\sigma }^{\text{fin,min}}}s_{h})^{\ell (f)}uu^{-1}
(\sum_{i\in \mathcal{P}_{\sigma }^{(2)}\cap \mathcal{P}_{\sigma }^{\text{fin,min}}}s_{i})^{\ast \ell (g)}\\
& = \delta _{\sigma }(f)u(\sum_{p\in \mathcal{P}_{\sigma }^{(\ell (f)+1)}\cap \mathcal{P}_{\sigma }^{\text{fin,max}}}s_{p})
(\sum_{q\in \mathcal{P}_{\sigma }^{(\ell (g)+1)}\cap \mathcal{P}_{\sigma }^{\text{fin,max}}}s_{q})u^{-1} & 
\end{align*}
and we conclude the proof.
\end{proof}

For any $f\in \mathcal{P}_{\sigma }^{\text{fin}}$,  
we let $\eta (f)$ denote the number $\text{min}\{\nu (z)\mid z\in V_{\sigma }^{(\ell (f)-1)}\}$. 
The following lemma holds.
\begin{lemma}
For any $f,g\in \mathcal{P}_{\sigma }^{\text{fin}}$ and $a,b\geq 0$, there exists $N\geq \ell (f)$ such that for any $n\geq N$, 
$u^{a}s_{f}s_{g}^{\ast }u^{-b}\in $ linear span $\{ u^{c}s_{p}s_{q}^{\ast }u^{-d} \mid 
p\in \mathcal{P}_{\sigma }^{(n)}$, $q\in  \mathcal{P}_{\sigma }^{(n-\ell (f)+\ell (g))}$, 
$0\leq c<\eta (p)$ and $0\leq d<\eta (q)\}$. 
\label{2}
\end{lemma}

\begin{proof}
According to the relation $u^{a}s_{f}s_{g}^{\ast }u^{-b}=\sum_{h\in \mathcal{P}_{\sigma }^{(n)}}u^{a}s_{f}s_{h}s_{h}^{\ast }s_{g}^{\ast }u^{-b}$ 
$(a,b\geq 0$, $n\geq 2$ and $f,g\in \mathcal{P}_{\sigma }^{\text{fin}})$ and   
$\lim_{n\rightarrow \infty }\text{min}\{\nu (z)\mid z\in V_{\sigma }^{(n)}\}=\infty $ 
 as $\sigma$ is primitive, the lemma holds.
\end{proof}

For a proper and primitive substitution $\sigma$, we let 
$$
\mathcal{A}_{\sigma }=\text{linear span} \{ u^{a}s_{f}s_{g}^{\ast }u^{-b} \mid a,b\geq 0 ~\text{and}~f,g\in \mathcal{P}_{\sigma }^{\text{fin}}\}
$$
and
$$
\mathcal{L}_{\sigma }=\{ u^{a}s_{f}s_{g}^{\ast }u^{-b} \mid 
f,g\in \mathcal{P}_{\sigma }^{\text{fin}}, 
0\leq a<\eta (f) ~\text{and}~ 0\leq b<\eta (g)\}.
$$ 

According to the Lemma \ref{2}, the following lemma holds.
\begin{lemma}
$$
\mathcal{A}_{\sigma }=\text{linear span} ~\mathcal{L}_{\sigma }.
$$
\end{lemma}

\begin{lemma}
Let $\sigma$ be a proper primitive substitution, then 
$$
\text{Alg}_{\mathbb{C}}\{u,\{s_{f} \}_{f\in E_{\sigma }^{(2)}} \}=\mathcal{A}_{\sigma }
$$
\label{lem:3}
\end{lemma}
\begin{proof}
Since $u,\{s_{f}\}_{f\in E_{\sigma }^{(2)}}$ are in $\mathcal{A}_{\sigma }$, 
to prove the lemma, 
it suffices to show that $\mathcal{A}_{\sigma }$ is an involutive algebra. 
We let $X=u^{a_{1}}s_{f_{1}}s_{g_{1}}^{\ast }u^{-b_{1}}, 
Y=u^{a_{2}}s_{f_{2}}s_{g_{2}}^{\ast }u^{-b_{2}}\in \mathcal{L}_{\sigma }$ with 
$f_{1},f_{2},g_{1}, g_{2}\in \mathcal{P}_{\sigma }^{\text{fin,min}}$. 
According to the Lemma \ref{2}, we may assume that $\ell (g_{1})=\ell (f_{2})$. We have
$$ 
XY=u^{a_{1}}s_{f_{1}}s_{g_{1}}^{\ast }u^{a_{2}-b_{1}}s_{f_{2}}s_{g_{2}}^{\ast }u^{-b_{2}}.
$$

We assume $a_{2}-b_{1}\geq 0$ without loss of generality, then 
$0\leq a_{2}-b_{1}<\nu(r(f_{2}))+ \eta (f_{2})$. 

If $a_{2}-b_{1}\leq \nu (r(f_{2}))-1$, then by using the Lemma \ref{1} (c), 
$XY
=u^{a_{1}}s_{f_{1}}s_{g_{1}}^{\ast }s_{\widehat{\lambda }_{\sigma }^{a_{2}-b_{1}}(f_{2})}s_{g_{2}}^{\ast }u^{-b_{2}}\in \mathcal{A}_{\sigma }$.

If $\nu (r(f_{2}))\leq a_{2}-b_{1}< \nu(r(f_{2}))+ \eta (f_{2})$, then 
we apply the Lemma \ref{1} (f) to $s_{f_{1}}s_{g_{1}}^{\ast }$, 
we have 
\begin{align*}
& XY=\sum_{h\in \mathcal{P}_{\sigma }^{(2)}\backslash \mathcal{P}_{\sigma }^{\text{fin,max}}}u^{a_{1}}s_{f_{1}}s_{\widehat{\lambda}_{\sigma } (h)}s_{\widehat{\lambda}_{\sigma } (h)}^{\ast }s_{g_{1}}^{\ast }u^{a_{2}-b_{1}}s_{f_{2}}s_{g_{2}}^{\ast }u^{-b_{2}}\\
& + \delta _{\sigma }(f_{1})u^{a_{1}+1}(\sum_{p\in \mathcal{P}_{\sigma }^{(\ell (f_{1})+1)}\cap \mathcal{P}_{\sigma }^{\text{fin,max}}}s_{p})
(\sum_{q\in \mathcal{P}_{\sigma }^{(\ell (g_{1})+1)}\cap \mathcal{P}_{\sigma }^{\text{fin,max}}}s_{q})^{\ast }u^{a_{2}-b_{1}-1}
s_{f_{2}}s_{g_{2}}^{\ast }u^{-b_{2}}.
\end{align*}

It follows from the Lemma \ref{1} (d) that 
$$
s_{\widehat{\lambda}_{\sigma } (h)}^{\ast }s_{g_{1}}^{\ast }u^{a_{2}-b_{1}}s_{f_{2}},~ 
s_{q}^{\ast }u^{a_{2}-b_{1}-1}s_{f_{2}}
\in \{s_{l}^{\ast }s_{r}\mid 
l,r\in \mathcal{P}_{\sigma }^{\text{fin}}\}
$$
for any $h\in \mathcal{P}_{\sigma }^{(2)}\backslash \mathcal{P}_{\sigma }^{\text{fin,max}}$ and 
$q\in \mathcal{P}_{\sigma }^{(\ell (g_{1})+1)}\cap \mathcal{P}_{\sigma }^{\text{fin,max}}$ and we conclude the proof.
\end{proof}

We let $D_{\sigma }$ denote the $C^{\ast }$-algebra generated by 
$\{u^{k}s_{f}s_{f}^{\ast }u^{-k}\mid f\in \mathcal{P}_{\sigma }^{\text{fin}}, 0\leq k\}.$

It is clear that 
any two of the form $s_{f}s_{f}^{\ast }$ $(f\in \mathcal{P}_{\sigma }^{\text{fin}})$ commute.
For $N\geq 2$, we let $D_{\sigma }^{(N)}$ denote the subalgebra of $D_{\sigma }$ generated by $\{s_{f}s_{f}^{\ast }\mid f\in \mathcal{P}_{\sigma }^{(N)}\}$. It is easily seen that $D_{\sigma }^{(N+1)}$ contains $D_{\sigma }^{(N)}$.

\begin{notation}
Let $\sigma$ be a proper primitive substitution. For any $f\in \mathcal{P}_{\sigma }^{\text{fin}}$, we let $\widetilde{U}(f)$ denote the set defined by
$$
\widetilde{U}(f)=\{g\in \mathcal{P}_{\sigma }^{\text{fin}}\mid 
g^{(1)}=f^{(1)},g^{(2)}=f^{(2)},\dots,g^{(\ell (f))}=f^{(\ell (f))}\}.
$$
\end{notation}

\begin{lemma}
For any $f\in \mathcal{P}_{\sigma }^{\text{fin}}$ and $a\geq 0$, there exists $N\geq \ell (f)$ such that for any $n\geq N$, 
$$
u^{a}s_{f}s_{f}^{\ast }u^{-a}\in  \text{linear span} \{ s_{g}s_{g}^{\ast }\mid 
g\in \mathcal{P}_{\sigma }^{(n)}\}.
$$
\label{5}
\end{lemma}
\begin{proof}
According to the Lemma \ref{1} (c), 
it suffices to show that for any $f\in \mathcal{P}_{\sigma }^{\text{fin,max}}$, 
there exists $n\in \mathbb{N}$ such that 
$us_{f}s_{f}^{\ast }u^{-1}\in $~linear span $\{ s_{g}s_{g}^{\ast }\mid g\in \mathcal{P}_{\sigma }^{(n)}\}$.
Since $\sigma $ is proper, we write $\mathcal{P}_{\sigma }^{\text{max}}=\{p_{\text{max}}\}$.
If $f\in \mathcal{P}_{\sigma }^{\text{fin,max}}$ with $p_{\text{max}}\not\in U(f)$, 
then we let $S(f)$ denote the set
$$
S(f)=\{g\in \widetilde{U}(f)\mid 
g^{(1)}, g^{(2)},\dots,g^{(\ell (g)-1)}\in E_{\sigma }^{\text{max}}~ \text{and}~ g^{(\ell (g))}\not\in E_{\sigma }^{\text{max}}\}.
$$

Then  by using the Cuntz-Krieger relations, we have
\begin{align*}
us_{f}s_{f}^{\ast }u^{-1}
 & = \sum_{g\in S(f)}us_{g}s_{g}^{\ast }u^{-1}
= \sum_{g\in S(f)}s_{\widehat{\lambda} _{\sigma }(g)}s_{\widehat{\lambda} _{\sigma }(g)}^{\ast }\\
 & \in \text{linear span}\{ s_{h}s_{h}^{\ast }\mid 
h\in \mathcal{P}_{\sigma }^{(N)}\} ~(\exists N\geq \ell (f)).
\end{align*}

If $f\in \mathcal{P}_{\sigma }^{\text{fin,max}}$ with $p_{\text{max}}\in U(f)$, then 
\begin{align*}
us_{f}s_{f}^{\ast }u^{-1}
& = u(1-\sum_{g\in \mathcal{P}_{\sigma }^{(\ell (f))},f\not=g}s_{g}s_{g}^{\ast })u^{-1}\\
& = 1-\sum_{g\in \mathcal{P}_{\sigma }^{(\ell (f))},f\not=g}s_{\widehat{\lambda _{\sigma }}(g)}s_{\widehat{\lambda _{\sigma }}(g)}^{\ast }
\end{align*}
and we conclude the proof.
\end{proof}

According to the Lemma \ref{5}, the following lemma holds.
\begin{lemma}
When $\sigma$ is proper and primitive,  
$D_{\sigma }$ is the inductive limit of the subalgebra $D_{\sigma }^{(N)}$.
\end{lemma}

\begin{lemma}
Let $\sigma$ be a proper primitive substitution. Then 
$C^{\ast }$-algebra $D_{\sigma }$ is isomorphic to $C(\mathcal{P}_{\sigma })$.
\label{lem:path}
\end{lemma}

\begin{proof}
The fact that $\mathcal{P}_{\sigma }$ is totally disconnected implies that linear span 
$\{ \chi _{f}\mid f\in \mathcal{P}_{\sigma }^{\text{fin}} \}$ is dense in 
$C(\mathcal{P}_{\sigma })$. 
It clearly holds that 
$\text{linear~span}\{\chi _{f}\mid f\in \mathcal{P}_{\sigma }^{(n)}\}$ is isomorphic to $D_{\sigma }^{(n)}$. 
Since $D_{\sigma }$ is the inductive limit of $D_{\sigma }^{(n)}$,  
$D_{\sigma }$ is isomorphic to $C(\mathcal{P}_{\sigma })$. 
\end{proof}

According to the Definition \ref{def:definition}, we have the following lemma.
\begin{lemma} 
Let $\sigma$ be a proper primitive substitution. 
In this lemma, we put $e_{f}=s_{f}s_{f}^{\ast }$, then 
the unitary $u$ and projections $\{e_{f}\}_{f\in \mathcal{P}_{\sigma }^{\text{fin}}}$ 
satisfy the following relations.\\
(i) $\sum_{f\in E_{\sigma }^{(2)}}e_{f}=1.$

For any $g\in \mathcal{P}_{\sigma }^{(n)}$, 
$e_{g}=\sum _{h\in E_{\sigma }^{(2)},~L(r(g))=L(s(h))}e_{gh}$.\\
(ii)
$$
e_{\widehat{\lambda }_{\sigma }(f)}=ue_{f}u^{-1}
$$
for any $n\geq 2$ and $f\in \mathcal{P}_{\sigma }^{(n)}\backslash 
(\mathcal{P}_{\sigma }^{(n)}\cap \mathcal{P}_{\sigma }^{\text{fin,max}})$.\\
(iii) 
$$
u(\sum_{f\in \mathcal{P}_{\sigma }^{(n)}\cap \mathcal{P}_{\sigma }^{\text{fin,max}}}e_{f})u^{-1}
=\sum_{f\in \mathcal{P}_{\sigma }^{(n)}\cap \mathcal{P}_{\sigma }^{\text{fin,min}}}e_{f}.
$$
for any $n\geq 2$.
\label{lem:4}
\end{lemma}

According to the Lemma \ref{lem:4}, we define the following algebra.
\begin{definition}
Let $\sigma$ be a proper primitive substitution. 
Then we define $\mathcal{F}_{\sigma }$ to be the universal unital $C^{\ast }$-algebra 
generated by a unitary $u$ and projections $\{e_{f}\}_{f\in \mathcal{P}_{\sigma }^{\text{fin}}}$ 
satisfying the relations (i), (ii) and (iii) in Lemma \ref{lem:4}.
\label{definition:def2}
\end{definition}

Let $(\mathcal{P}_{\sigma }, \lambda _{\sigma })$ be a Bratteli-Vershik system 
associated with a proper primitive substitution $\sigma$.  
Then this gives a $C^{\ast }$-dynamical system 
$(C(\mathcal{P}_{\sigma })$, $\mathbb{Z}$, $\widetilde{\lambda} _{\sigma })$, where 
$$
\widetilde{\lambda} _{\sigma }^{n}(f)=f\circ \lambda _{\sigma }^{-n} ~(n\in \mathbb{Z}).
$$
\begin{prop} 
Let $\sigma$ be a proper primitive substitution. Then $\mathcal{F}_{\sigma }$ is 
isomorphic to the crossed product  
$C(\mathcal{P}_{\sigma })\rtimes _{\widetilde{\lambda} _{\sigma }}\mathbb{Z}$.
\end{prop}
\begin{proof}
The crossed product is generated by projections $\chi _{f}$ $(f\in \mathcal{P}_{\sigma }^{\text{fin}})$ and a unitary $U$ implementing the automorphism $\widetilde{\lambda} _{\sigma }$.
It follows from Lemma \ref{lem:path} that 
there is a homomorphism of $C(\mathcal{P}_{\sigma })$ into $\mathcal{F}_{\sigma }$ such that $\phi (\chi _{f})=e_{f}$ for any $f\in \mathcal{P}_{\sigma }^{\text{fin}}$. 
We check that 
$$
\phi (\widetilde{\lambda} _{\sigma }(\chi _{f}))=ue_{f}u^{-1}
$$
for any $f\in \mathcal{P}_{\sigma }^{\text{fin}}$. 
If $f\in \mathcal{P}_{\sigma }^{\text{fin}}\backslash 
\mathcal{P}_{\sigma }^{\text{fin,max}}$, then $\phi (\widetilde{\lambda }_{\sigma }(\chi _{f}))
=e_{\widehat{\lambda }_{\sigma }(f)}=ue_{f}u^{-1}.$ 
We write $\mathcal{P}_{\sigma }^{\text{max}}=\{p_{\text{max}}\}$. 
If $f\in \mathcal{P}_{\sigma }^{\text{fin,max}}$ is such that $p_{\text{max}}\not\in U(f)$, 
then we take the set $S(f)\subset \mathcal{P}_{\sigma }^{\text{fin}}$ (The definition of $S(f)$ is in the proof of Lemma \ref{5}.). 
Then $\chi _{f}=\sum_{g\in S(f)}\chi _{g}$. So
$$
\phi (\widetilde{\lambda }_{\sigma }(\chi _{f}))
 = \sum_{g\in S(f)}e_{\widehat{\lambda }_{\sigma }(g)}
 = u(\sum_{g\in S(f)}e_{g})u^{-1}=ue_{f}u^{-1}.
$$

If $f\in \mathcal{P}_{\sigma }^{\text{fin,max}}$ such that $p_{\text{max}}\in U(f)$, then $\chi _{f}=1-\sum_{g\in \mathcal{P}_{\sigma }^{(\ell (f))},f\not=g}\chi _{g}$. So 
$$
\phi (\widetilde{\lambda} _{\sigma }(\chi _{f}))
 = \phi (\widetilde{\lambda} _{\sigma }(1-\sum_{g\in \mathcal{P}_{\sigma }^{(\ell (f))},f\not=g}\chi _{g}))
 =1-\sum_{g\in \mathcal{P}_{\sigma }^{(\ell (f))},f\not=g}ue_{g}u^{-1}= ue_{f}u^{-1}.
$$

It follows from the discussion above 
that $\mathcal{F}_{\sigma }$ provides a covariant representation of $(C(\mathcal{P}_{\sigma }), \mathbb{Z}, \widetilde{\lambda} _{\sigma })$.
It follows from the universal property of the crossed product 
there is a homomorphism of $C(\mathcal{P}_{\sigma })\rtimes _{\widetilde{\lambda} _{\sigma }}\mathbb{Z}$ 
onto $\mathcal{F}_{\sigma }$ taking $U$ to $u$ and $\chi _{f}$ to $e_{f}$. Since 
$\lambda _{\sigma }$ is minimal, 
$C(\mathcal{P}_{\sigma })\rtimes _{\widetilde{\lambda }_{\sigma }}\mathbb{Z}$ is simple and 
we conclude the proof.
\end{proof}

We consider a unitary $u$ and partial isometries $\{s_{f}\}_{f\in E_{\sigma }^{(2)}}$ 
satisfying the relations (a)-(c) of 
Definition \ref{def:definition} which generate $B_{\sigma }$. 
By the universal property of $B_{\sigma }$, for any $t\in \mathbb{T}$, 
there exists an automorphism $\alpha _{t}$ of $B_{\sigma }$ such that 
$\alpha _{t}(u)=u$, $\alpha _{t}(s_{f})=ts_{f}$ for each 
$f$ in $E_{\sigma }^{(2)}$.
 Put
$$
E(x)=\int _{\mathbb{T}}\alpha _{t}(x)dt,~~~~~~x\in B_{\sigma }
$$
where $dt$ is the normalized Haar measure on $\mathbb{T}$. 
Then $E$ is a conditional expectation onto $\mathcal{F}_{\sigma }$.

We consider a unitary $u$ and partial isometries $\{e_{f}\}_{f\in E_{\sigma }^{(2)}}$ 
satisfying the relations (i)-(iii) of 
Lemma \ref{lem:4} which generate $\mathcal{F}_{\sigma }$. 
By the universal property of $B_{\sigma }$, for any $t\in \mathbb{T}$, 
there exists an automorphism $\beta _{t}$ of $B_{\sigma }$ such that 
$\beta _{t}(u)=e^{it}u$, $\beta _{t}(e_{f})=e_{f}$ for each 
$f$ in $E_{\sigma }^{(2)}$. Put
$$
F(x)=\int_{\mathbb{T}}\beta _{t}(x)dt,~~~~~x\in \mathcal{F}_{\sigma }.
$$

Then $F$ is a conditional expectation onto $D_{\sigma }$.

We put $G=F\circ E$, then $G$ is a conditional expectation of $B_{\sigma }$ onto  $D_{\sigma }$.

\begin{thm} Let $\sigma$ be a proper primitive substitution. Then $C^{\ast }$-
algebra $B_{\sigma }$ associated with $\sigma$ is simple and purely infinite.
\label{th:simple}
\end{thm}
\begin{proof}
Let $\sigma$ be a proper primitive substitution on an alphabet $A$ with $^\# A=c$.
To prove the theorem, it suffices to show that for any non-zero positive element 
$X\in B_{\sigma }$, there exists $Z\in B_{\sigma }$ such that 
$Z^{\ast }XZ$ is invertible.
Since $G$ is faithful, $G(X)\not=0$. We may assume that $\| G(X)\| =1$.
For arbitrary $0<\varepsilon <\frac{1}{3c^{2}+1}$, there exists 
$Y\in $
$\text{Alg}_{\mathbb{C}}$
$\{u,\{s_{f}\}_{f\in E_{\sigma }^{(2)}} \}$
such that $\|  X-Y \| <\varepsilon $. Then we have
$$
\| G(Y)\| =\| G(X)-G(X-Y)\| \geq \| G(X)\|-\| X-Y\|  >\frac{3c^{2}}{3c^{2}+1}.
$$

For any $N\geq 2$, we let 
\begin{align*}
& \mathcal{L}_{\sigma }^{(N)}=\{u^{a}s_{f}s_{g}^{\ast }u^{-b}\mid f,g\in \mathcal{P}_{\sigma }^{\text{fin,min}}, 
\text{min}\{\ell (f), \ell (g)\}=N, \\
& \text{max}\{\ell (f), 
\ell (g)\}\leq 2N  ~\text{and}~ \text{if}~ u^{a}s_{f}s_{g}^{\ast }u^{-b}\in D_{\sigma }, \text{then}~ 
 u^{a}s_{f}s_{g}^{\ast }u^{-b}\in D_{\sigma }^{(N)}.\}
\end{align*}
and $\mathcal{A}_{\sigma }^{(N)}=\text{linear span}~\mathcal{L}_{\sigma }^{(N)}$.
According to the Lemma \ref{2} and Lemma \ref{5}, 
there exists $N>0$ such that $Y\in \mathcal{A}_{\sigma }^{(N)}$. 
Now we choose pathes $\{p(f,v)|f\in \mathcal{P}_{\sigma }^{(N)}, 
v\in V_{\sigma }^{(1)}\}$ $\subset \mathcal{P}_{\sigma }^{\text{fin}}$ as follows:
$$
p(f,v)=p_{1}(f,v)p_{2}(f,v)\in \mathcal{P}_{\sigma }^{\text{fin}}
$$
where 
$p_{1}(f,v)\in \widetilde{U}(f)$ is such that 
when we write $p_{1}(f,v)=\widehat{\lambda} _{\sigma }^{k}(g) ~(g\in \mathcal{P}_{\sigma }^{\text{fin,min}} ~\text{with}~r(g)=r(p_{1}(f,v))~\text{and}~
0\leq k  )$, $\text{max}\{\nu (z)\mid z\in V_{\sigma }^{(2N)}\}\leq k\}$.

$p_{2}(f,v)\in \mathcal{P}_{\sigma }^{\text{fin}}$ is a sufficiently large length path such that $L(s(p_{2}(f,v)))=L(r(p_{1}(f,v)))$, 
$L(r(p_{2}(f,v)))=L(v)$ and 
$s_{p_{2}(f,v)}^{\ast }s_{g}s_{p_{2}(f,v)}=0$ 
for any $g\in \cup _{m=2}^{N+1}\mathcal{P}_{\sigma }^{(m)}$.

We claim that if $x=u^{a}s_{g}s_{h}^{\ast }u^{-b}\in \mathcal{L}^{(N)}_{\sigma }$ 
$(g,h\in \mathcal{P}_{\sigma }^{\text{fin,min}})$ 
with 
$G(x)=0$, then $e_{p(f,v)}xe_{p(f,v)}=0$ for any $f\in \mathcal{P}_{\sigma }^{(N)}$ and  $v\in V_{\sigma }^{(1)}$, because
\begin{align*}
& e_{p(f,v)}xe_{p(f,v)} = s_{p(f,v)}s_{\widehat{\lambda }_{\sigma }^{-a}(p(f,v))}^{\ast }s_{g}s_{h}^{\ast }
s_{\widehat{\lambda }_{\sigma }^{-b}(p(f,v))}s_{p(f,v)}^{\ast }\\
& = s_{p(f,v)}s_{p_{2}(f,v)}^{\ast }s_{\widehat{\lambda }_{\sigma }^{-a}(p_{1}(f,v))}^{\ast }s_{g}s_{h}^{\ast }
s_{\widehat{\lambda }_{\sigma }^{-b}(p_{1}(f,v))}s_{p_{2}(f,v)}s_{p(f,v)}^{\ast }=0.
\end{align*}

For any $v\in V_{\sigma }^{(1)}$, 
we let $\varphi _{v}^{(N)}$: $D_{\sigma }\longrightarrow D_{\sigma }$ be a map defined by
\begin{equation*}
\varphi _{v}^{(N)}(x)=\sum_{f\in \mathcal{P}_{\sigma }^{(N)}}e_{p(f,v)}xe_{p(f,v)}.
\end{equation*}

Then the restlicted map 
$$
\varphi _{v}^{(N)}|_{D_{\sigma }^{(N)}}:D_{\sigma }^{(N)}\longrightarrow 
C^{\ast }(\{e_{p(f,v)}|f\in \mathcal{P}_{\sigma }^{(N)}\})
$$
is an isomorphism. 
According to the discussion above and the definition of $\mathcal{L}_{\sigma }^{(N)}$, we have
$$
\varphi _{v}^{(N)}(G(Y))=\sum_{f\in \mathcal{P}_{\sigma }^{(N)}}e_{p(f,v)}Ye_{p(f,v)}.
$$

So there exists a path 
$p_{v}\in \{p(f,v)|f\in \mathcal{P}_{\sigma }^{(N)}\}$
such that
\begin{equation*}
e_{p_{v}}Ye_{p_{v}}=\| G(Y)\|e_{p_{v}}
\end{equation*}
for each $v\in V_{\sigma }^{(1)}$. 

Put $Z=\sum_{v\in V_{\sigma }^{(1)}}\| G(Y)\| ^{-\frac{1}{2}}$
$s_{p_{v}}$, then $\| Z\| \leq \| G(Y)\|^{-\frac{1}{2}}c$ and we have
$$
Z^{\ast }YZ=\sum_{v\in V_{\sigma }^{(1)}}s_{p_{{v}}}^{\ast}s_{p_{v}}=1.
$$

Then 
$\| 1-Z^{\ast }XZ\| =\|Z^{\ast }YZ-Z^{\ast }XZ \| \leq \| Z\|^{2}\| Y-X\|$
$<\frac{1}{3}$.

Thus $Z^{\ast }XZ$ is invertible and we conclude the proof.
\end{proof}

\begin{remark}
Let $\sigma$ be a proper primitive substitution and $n\in \mathbb{N}$. 
Then $\{s_{f}\}_{f\in \mathcal{P}_{\sigma }^{(n+1)}}\subset B_{\sigma }$ satisfy the  relations $(a)$-$(c)$ in Definition \ref{def:definition} of $B_{\sigma ^{n}}$. 
Thus $B_{\sigma ^{n}}$ is isomorphic to $B_{\sigma }$.
\label{remark}
\end{remark}
\section{actions for $B_{\sigma }$}
Let $\sigma$ be a primitive substitution on an alphabet $A$. 
For $x\in A^{\mathbb{Z}}$, we denote by $x_{[j,~k]}$ the word $x_{j}x_{j+1}\dots x_{k}$.
We define
\begin{align*}
X_{\sigma } =
 &  \{(x_{i})\in A^{\mathbb{Z}}\mid \text{For any} ~j<k, ~\text{there exists}~ n\in \mathbb{N} ~\text{and} ~a\in A\\
& \text{such that}~ x_{[j,k]}~ \text{is a subword of}~ \sigma ^{n}(a).\}
\end{align*}
and let $T_{\sigma }:X_{\sigma }\longrightarrow X_{\sigma }$ be a map defined by 
$T_{\sigma }((x_{k})_{k\in \mathbb{Z}})=(x_{k+1})_{k\in \mathbb{Z}}$.
It is known in \cite{Qu} that for any primitive substitution $\sigma $, the substitution dynamical system $(X_{\sigma },T_{\sigma })$ is minimal and uniquely ergodic.

A substitution $\sigma :A\longrightarrow \{\text{the set of words of}$ $A \}$ naturally extends to a map 
$\sigma :A^{\mathbb{Z}}\longrightarrow A^{\mathbb{Z}}$:
$$
\sigma (...w_{-2}w_{-1}.w_{0}w_{1}...)=...\sigma (w_{-2})\sigma (w_{-1}).\sigma (w_{0})\sigma (w_{1})....
$$

Let $\sigma$ be a proper primitive substitution.  
An infinite word $w\in A^{\mathbb{Z}}$ is  called a periodic point of $\sigma$ if $\sigma (w)=w$ and is $T_{\sigma }$-periodic if there exists a positive integer $k$ with 
$T_{\sigma }^{k}(w)=w$. 
A proper substitution $\sigma$ is said to be periodic when there exists a periodic point of $\sigma$ when is also $T_{\sigma }$-periodic. When $\sigma$ is primitive, $|X_{\sigma }|<\infty$ if and only if $\sigma$ is periodic.

If $\sigma$ is periodic proper and primitive,  then we denote by $p_{\sigma }$ the smallest period of its fixed point $w$ and $u=w_{[0,~p_{\sigma }-1]}$. 
As $\sigma (w)=w$, there exists $d_{\sigma }\geq 1$ such that 
$\sigma (u)=u\dots u$ $(d_{\sigma }~times)$.  
We need the following fact due to \cite{Du}.~$(\text{Proposition 16})$.
\begin{prop}
Let $\sigma$ be a proper primitive substitution.\\  
(a)If $\sigma$ is non-periodic, then the Bratteli-Vershik 
system $(\mathcal{P}_{\sigma },\lambda _{\sigma }) $ is isomorphic to the system $(X_{\sigma },T_{\sigma })$.\\
(b)If $\sigma $ is periodic, the system $(\mathcal{P}_{\sigma },\lambda _{\sigma })$ is 
isomorphic to an odometer with a stationary base $(p_{\sigma },d_{\sigma },d_{\sigma },\dots)$.
\end{prop}
The next lemma follows from \cite{Da}$(\text{Theorem~VIII.4.1})$.
\begin{lemma}
If $\sigma$ is a periodic, proper  and primitive substitution, then 
$\mathcal{F}_{\sigma }$ is isomorphic to Bunce-Deddens algebra with supernatural number
$p_{\sigma }d_{\sigma }^{\infty }$.
\end{lemma}

\begin{definition}
Let $\sigma$ be a primitive substitution on an alphabet $A$. 
For any $a, b\in A$, let $\kappa (a,b)\in \mathbb{N}$ denote the number of occurrences of $b$ in $\sigma (a)$. 
Let $M_{\sigma }$ denote its matrix; $M_{\sigma }=\{\kappa (a,b)\}_{a,b\in A}$. 
A primitive substitution $\sigma$ is para-periodic if there exists $p\in \mathbb{N}$ such that $M_{\sigma }^{2}=pM_{\sigma }$. Let $\sigma$ be a para-periodic substitution, then we define a number $N_{\sigma }\in \mathbb{N}$ by 
$M_{\sigma }^{2}=N_{\sigma }M_{\sigma }.$
\end{definition}

Let $\sigma$ be a proper primitive substitution. 
It is known in \cite{Da} $(\text{Corollary~VIII.3.8})$ that there is a unique faithful trace $\varphi _{\sigma }^{(0)}$ 
on $\mathcal{F}_{\sigma }$.
 We let $\varphi_{\sigma }$ be a state on $B_{\sigma }$ defined by 
$\varphi _{\sigma }=\varphi _{\sigma }^{(0)}\circ E$.

\begin{lemma}
Let $\sigma$ be a proper primitive substitution on an alphabet $A$.\\ 
(a)For any $f,g,h\in \mathcal{P}_{\sigma }^{\text{fin}}$ and $0\not=k\in \mathbb{Z}$,
$$
\varphi _{\sigma }(u^{k}s_{f}s_{f}^{\ast })=
\varphi _{\sigma }(s_{g}s_{h}^{\ast }u^{k}s_{h}s_{g}^{\ast })=0.
$$
(b)If $\sigma$ is para-periodic, then  
for any $f,g\in \mathcal{P}_{\sigma }^{\text{fin}}$ with $L(r(f))=L(r(g))$, we have

$$
\varphi _{\sigma }(s_{f}s_{f}^{\ast })=N_{\sigma }^{\ell (g)-\ell (f)}
\varphi _{\sigma }(s_{g}s_{g}^{\ast }).
$$
\label{7}
\end{lemma}

\begin{proof}
(a)Let $f,g,h\in \mathcal{P}_{\sigma }^{\text{fin}}$. 
According to the relation 
$$
u^{k}s_{f}s_{f}^{\ast }=\sum_{g\in \mathcal{P}_{\sigma }^{(n)}}
u^{k}s_{f}s_{g}s_{g}^{\ast }s_{f}^{\ast }
$$ 
for any $n\geq 2$, 
we may first assume that $|k|<\nu (r(f))$ and $k\not=0$. Then $s_{f}^{\ast }u^{k}s_{f}=0$.
Thus
$$
\varphi _{\sigma }(u^{k}s_{f}s_{f}^{\ast })=
\varphi _{\sigma }(s_{f}s_{f}^{\ast }u^{k}s_{f}s_{f}^{\ast })=0.
$$ 
$$
\varphi _{\sigma }(s_{g}s_{h}^{\ast }u^{k}s_{h}s_{g}^{\ast })
=\sum_{p\in \mathcal{P}_{\sigma }^{(n)}}
\varphi _{\sigma }(s_{p}s_{p}^{\ast }s_{g}s_{h}^{\ast }u^{k}s_{h}s_{g}^{\ast }
s_{p}s_{p}^{\ast })=0~(n\geq 2).
$$
(b) 
We let $p_{b}^{(n)}$ denote the path in $\mathcal{P}_{\sigma }^{(n)}\cap \mathcal{P}_{\sigma }^{\text{fin,min}}$ such that $L(r(p_{b}^{(n)}))=b$ for any $b\in A$ 
and $n\in \mathbb{N}$. 
Let $f,g\in \mathcal{P}_{\sigma }^{\text{fin}}$ with $L(r(f))=L(r(g))$ and $\ell (f)\leq \ell (g)$, then
$$
\varphi _{\sigma }(s_{f}s_{f}^{\ast })=\sum_{h\in \mathcal{P}_{\sigma }^{(\ell (g)-\ell (f)+2)}}\varphi _{\sigma }(s_{f}s_{h}s_{h}^{\ast }s_{f}^{\ast }).
$$

Since 
$\varphi _{\sigma }(s_{f}s_{f}^{\ast })=\varphi _{\sigma }(s_{f^{'}}s_{f^{'}}^{\ast })$ 
for any $f,f^{'}\in \mathcal{P}_{\sigma }^{\text{fin}}$ with $r(f)=r(f^{'})$ and 
\begin{align*}
\ ^\#\{h\in \mathcal{P}^{(\ell (g)-\ell (f)+2)}\mid L(s(h))=L(r(f)) ~\text{and} ~L(r(h))=b\}=N_{\sigma }^{\ell (g)-\ell (f)}\kappa (b,L(r(f)))
\end{align*}
as $\sigma$ is para-periodic, we have
\begin{align*}
\sum_{h\in \mathcal{P}_{\sigma }^{(\ell (g)-\ell (f)+2)}}\varphi _{\sigma }(s_{f}s_{h}s_{h}^{\ast }s_{f}^{\ast })
& = N_{\sigma }^{\ell (g)-\ell (f)}
\sum_{b\in A}\kappa (b,L(r(f)))\varphi _{\sigma }(s_{p_{b}^{(\ell (g)+1)}}
s_{p_{b}^{(\ell (g)+1)}}^{\ast })\\
& = N_{\sigma }^{\ell (g)-\ell (f)}
\varphi _{\sigma }(s_{g}s_{g}^{\ast }).
\end{align*}
\end{proof}

\begin{lemma}
Let $\sigma$ be a proper primitive substitution. 
Let $x=u^{a}s_{f}s_{g}^{\ast }u^{-b}\in \mathcal{L}_{\sigma }$, where $f, g\in \mathcal{P}_{\sigma }^{\text{fin}, \text{min}}$. 
For any $n>\ell (f)$, $x$ is written as the form
$$ 
x=\sum_{\text{finite}}u^{a_{k}}s_{f_{k}}s_{g_{k}}^{\ast }u^{-b_{k}}
$$
where $u^{a_{k}}s_{f_{k}}s_{g_{k}}^{\ast }u^{-b_{k}}\in \mathcal{L}_{\sigma }$, 
$f_{k}\in \mathcal{P}_{\sigma }^{(n)}\cap 
\mathcal{P}_{\sigma }^{\text{fin,min}}$ $~\text{and}~ 
\eta (f_{k})\leq a_{k}< \nu (r(f_{k}))+\eta (f_{k})$.
 
\label{6}
\end{lemma}
\begin{proof}
Let $x=u^{a}s_{f}s_{g}^{\ast }u^{-b}\in \mathcal{L}_{\sigma }$ and $n>\ell (f)$. 
According to the relation $\sum_{h\in \mathcal{P}_{\sigma }^{(m)}}s_{h}s_{h}^{\ast }=1$
 $(m\geq 1)$, 
we may assume that $\ell (f)=n-1$. 
If $a\geq \text{min}\{\nu (z)\mid \ell (z)=n-1\}$, then 
$$
x=u^{a}s_{f}s_{g}^{\ast }u^{-b}=\sum_{h\in E_{\sigma }^{(2)}}u^{a}s_{f}s_{h}s_{h}^{\ast }s_{g}^{\ast }u^{-b}. 
$$

If $a<\text{min}\{\nu (z)\mid \ell (z)=n-1\}$, according to the Lemma 3.2 (f), we have
\begin{align*}
x
& = \sum_{h\in E_{\sigma }^{(2)}\backslash E_{\sigma }^{\text{max}}}u^{a}s_{f}
s_{\widehat{\lambda}_{\sigma } (h)}s_{\widehat{\lambda}_{\sigma } (h)}^{\ast }s_{g}^{\ast }u^{-b}\\
& + u^{a+1}(\sum_{h\in \mathcal{P}_{\sigma }^{(\ell (f)+1)}\cap \mathcal{P}_{\sigma }^{\text{fin,max}}}s_{h})
(\sum_{i\in \mathcal{P}_{\sigma }^{(\ell (g)+1)}\cap \mathcal{P}_{\sigma }^{\text{fin,max}}}s_{i})^{\ast }u^{-(b+1)}.
\end{align*}

The above equations satisfy the condition in Lemma \ref{6}.
\end{proof}

\begin{definition} Let $\mathcal{A}$ be a $C^{\ast }$-algebra, $\theta $ an $\mathbb{R}$-action on $\mathcal{A}$, $\varphi$ a state on $\mathcal{A}$ and $\beta>0$.
Then $\varphi$ satisfies the $\beta$-KMS condition if 
for each $x$, $y\in \mathcal{A}$, 
there is a function $f$, 
bounded and continuous on $\{z\in \mathbb{C}\mid 0\leq \text{Im}(z)\leq \beta \}$ and 
analytic on $\{z\in \mathbb{C}\mid 0<\text{Im}(z)<\beta \}$  such that for $t\in \mathbb{R}$,
$$
f(t)=\varphi (\theta _{t}(x)y),~~ f(t+i\beta )=\varphi (y\theta _{t}(x)).
$$
\end{definition}

It is known that $\varphi$ satisfies the $\beta$-KMS condition if and only if 
$\varphi (x\theta _{i\beta }(y))=\varphi (yx)$ for any $x\in \mathcal{A}$ and $y$, 
analytic element for $(\theta _{t})_{t\in \mathbb{R}}$.

Let $\sigma$ be a para-periodic proper substitution. Then let $\theta  :\mathbb{R}\longrightarrow \text{Aut}(B_{\sigma })$ be a map defined by
$$
\theta _{t}(u)=u,~~
\theta _{t}(s_{f})=N_{\sigma }^{it}s_{f}~~~(f\in E_{\sigma }^{(2)}).
$$

Then $\theta $ is an $\mathbb{R}$-action on $B_{\sigma }$.
\begin{prop} Let $\sigma$ be a para-periodic proper substitution. Then the state $\varphi_{\sigma }$ is a $1$-KMS state for $(\theta _{t})_{t\in \mathbb{R}}$.
\end{prop}
\begin{proof}
To prove the proposition, by Lemma \ref{lem:3} and Lemma \ref{6}, it suffices to show that 
$\varphi _{\sigma }(xy)=\varphi _{\sigma }(y\theta _{i}(x))$ for any 
$x=u^{a_{1}}s_{f_{1}}s_{g_{1}}^{\ast }u^{-b_{1}}$, 
$y=u^{a_{2}}s_{f_{2}}s_{g_{2}}^{\ast }u^{-b_{2}}\in \mathcal{L}_{\sigma }$ with 
$f_{1}$,$f_{2}$,$g_{1}$,$g_{2}\in \mathcal{P}_{\sigma }^{\text{fin,min}}$, $\ell (g_{1})=\ell (f_{2})$, $\eta (g_{1})\leq b_{1}< \nu (r(g_{1}))+\eta (g_{1})$ and 
$\eta (f_{2})\leq a_{2}< \nu (r(f_{2}))+\eta (f_{2})$. 

If $\ell (f_{1})+\ell (f_{2})\not=\ell (g_{1})+\ell (g_{2})$, then 
$xy,yx\not\in \mathcal{F}_{\sigma }$ and so $\varphi _{\sigma } (xy)=\varphi _{\sigma } (yx)=0$.

We let $\ell (f_{1})+\ell (f_{2})=\ell (g_{1})+\ell (g_{2})$, then 
$\ell (f_{1})=\ell (g_{2})$.
If $b_{1}\not=a_{2}$ or $f_{2}\not=g_{1}$, then $s_{g_{1}}^{\ast }u^{a_{2}-b_{1}}s_{f_{2}}=0$. Thus $xy=0$ and 
$\varphi _{\sigma } (yx)=\varphi _{\sigma } (s_{g_{1}}s_{g_{1}}^{\ast }u^{a_{2}-b_{1}}s_{f_{2}}s_{g_{2}}^{\ast }u^{a_{1}-b_{2}}s_{f_{1}}s_{g_{1}}^{\ast })=0$.

Next we let $b_{1}=a_{2}$ and $g_{1}=f_{2}$, then $f_{1}=g_{2}$. 
If $a_{1}\not=b_{2}$, then by using the Lemma \ref{7} (a), $\varphi _{\sigma } (xy)=\varphi _{\sigma }(u^{a_{1}-b_{2}}s_{f_{1}}s_{f_{1}}^{\ast })=0$ and $\varphi _{\sigma }(yx)=
\varphi _{\sigma }(s_{f_{1}}s_{g_{1}}^{\ast }u^{a_{1}-b_{2}}s_{g_{1}}s_{f_{1}}^{\ast })=0$.

We let $a_{1}=b_{2}$. According to the Lemma \ref{7} (b), we have 
\begin{align*}
\varphi _{\sigma }(xy)
& = \varphi _{\sigma }(s_{f_{1}}s_{f_{1}}^{\ast })\\
& = N_{\sigma }^{\ell (g_{1})-\ell (f_{1})}\varphi _{\sigma }(s_{g_{1}}s_{g_{1}}^{\ast })\\
& = \varphi _{\sigma }(y\theta _{i}(x))
\end{align*}
and we conclude the proof.
\end{proof}

\section{The $K$-groups of $B_{\sigma }$}
We follow the definition of semigroup crossed product $C^{\ast }$-algebras 
in~\cite{R}. 
\begin{definition} 
Let $B$ be a unital $C^{\ast }$-algebra and $\rho$ be an endomorphism of $B$ such that 
$\rho$ is an isomorphism from $B$ onto $\rho (1)B\rho (1)$.
We define $B\rtimes _{\rho }\mathbb{N}$ to be the universal unital 
$C^{\ast }$-algebra generated by $B$ and an isometry $s$ such that $\rho (b)=sbs^{\ast }$
for any $b\in B$.
\end{definition}
We denote by $\mu _{1}$ an isometry $\sum_{f\in E_{\sigma }^{(2)}\cap E_{\sigma }^{\text{min}}}s_{f}$ and let $\rho$ be a map defined by
\begin{equation*}
\rho :\mathcal{F}_{\sigma }\longrightarrow \mathcal{F}_{\sigma },~~~
\rho (x)=\mu _{1}x\mu _{1}^{\ast }.
\end{equation*}
Then $\rho \in \text{End}(\mathcal{F}_{\sigma })$ and it is easily seen that 
$B_{\sigma }$ is isomorphic to $\mathcal{F}_{\sigma }\rtimes _{\rho }\mathbb{N}$.

We shall recall some facts on this dimension group $K_{\sigma }$ of a substitution dynamical system based on \cite{Du}.

\begin{definition}
Let $\sigma$ be a primitive substitution on an alphabet $A$ with $^\#A=d$,  then we let\\ $G_{\sigma }=\{x\in \mathbb{Q}^{d}\mid \text{there exists}~n>0 ~\text{such that}
~M_{\sigma }^{n}x\in \mathbb{Z}^{d}.\}$.\\
$H_{\sigma }=\{x\in \mathbb{Q}^{d}\mid \text{there exists}~n>0 ~\text{such that}~
M_{\sigma }^{n}x=0. \}$.\\
$K_{\sigma }=G_{\sigma }/H_{\sigma }$.

The group $K_{\sigma }$ is called the dimension group of a substitution dynamical system. 
\end{definition}

We let $\phi$ be an endomorphism on $\mathbb{Z}^{d}$ defined by $\phi (x)=M_{\sigma }x$,  then the direct limit of $(\mathbb{Z}^{d},\phi )$ is isomorphic to $K_{\sigma }$.

Let $\sigma$ be a proper primitive substitution on an alphabet $A=\{c_{1},c_{2} \dots c_{d}\}$, then we let $C(\mathcal{P}_{\sigma },\mathbb{Z})$ denote the set of continuous functions with values in $\mathbb{Z}$. It follows from the Pimsner-Voiculescu exact sequence that 
$K_{0}(\mathcal{F}_{\sigma })=C(\mathcal{P}_{\sigma },\mathbb{Z})/\widetilde{\lambda} _{\sigma }C(\mathcal{P}_{\sigma },\mathbb{Z})$ 
where $\widetilde{\lambda} _{\sigma }(f)=f-f\circ \lambda _{\sigma }^{-1}$, for each $f\in C(\mathcal{P}_{\sigma },\mathbb{Z})$ and 
$K_{1}(\mathcal{F}_{\sigma })=\mathbb{Z}$.

For $n\geq 2$, we denote by $C_{n}$ the subgroup of $C(\mathcal{P}_{\sigma },\mathbb{Z})$ and $H_{n}$ the subgroup of $C_{n}$ defined by 
$$
C_{n}=\{\sum_{f\in \mathcal{P}_{\sigma }^{(n)}}a_{f}\chi _{f}\mid a_{f}\in \mathbb{Z}\}, 
$$
$$
H_{n}=\{\sum_{f\in \mathcal{P}_{\sigma }^{(n)}}a_{f}\chi _{f}\in C_{n}\mid 
 \sum_{f\in \mathcal{P}_{\sigma }^{(n)},~L(r(f))=c_{k}}a_{f}=0 ~\text{for any}~ 
1\leq k\leq d\}.
$$

Then it clearly holds that $C_{n}\subset C_{n+1}$, $H_{n}\subset H_{n+1}$, 
$C(\mathcal{P}_{\sigma },\mathbb{Z})=\cup _{n=2}^{\infty }C_{n}~\text{and}$ 
$\widetilde{\lambda} _{\sigma }C(\mathcal{P}_{\sigma },\mathbb{Z})=\cup _{n=2}^{\infty }H_{n}$. 
We let $K_{n}=C_{n}/H_{n}$ and 
$\varphi _{n}$ be a map on $K_{n}$ defined by 
\begin{align*}
\varphi _{n}:K_{n}\longrightarrow K_{n+1},~\varphi _{n}(x+H_{n})=x+H_{n+1}
\end{align*}
for $n\geq 2$. 

The group $C_{n}/H_{n}$ is isomorphic to $\mathbb{Z}^{d}$ for any $n\geq 2$ with 
isomorphism $\psi _{n}$ defined by 
\begin{align*}
\psi _{n}(\sum_{f\in \mathcal{P}_{\sigma }^{(n)}}a_{f}\chi _{f}+H_{n})=
(\sum_{g\in \mathcal{P}_{\sigma }^{(n)},L(r(g))=c_{k} }a_{g})_{1\leq k\leq d}.
\end{align*}

Then we have 
$\phi \circ \psi _{n}=\psi _{n+1} \circ \varphi _{n}$ on $K_{n}$~$(n\geq 2)$.
Since $C(\mathcal{P}_{\sigma },\mathbb{Z})/\widetilde{\lambda} _{\sigma }C(\mathcal{P}_{\sigma },\mathbb{Z})$ 
is the direct limit of $(C_{n}/H_{n},\phi _{n})$ and 
$K_{\sigma }$ is the direct limit of $(\mathbb{Z}^{d},\phi )$
 respectively, $K_{0}(\mathcal{F}_{\sigma })$ is isomorphic to $K_{\sigma }$.

Let $\sigma$ be a proper primitive substitution. For any $x+H_{\sigma }\in K_{\sigma }$, 
there exists 
$y+H_{\sigma }\in K_{\sigma }$ such that $x-M_{\sigma }y\in H_{\sigma }$. 
Then we let $\tau _{\sigma }:K_{\sigma }\longrightarrow K_{\sigma }$ be a map defined by 
$\tau _{\sigma }(x+H_{\sigma })=y+H_{\sigma }$.

\begin{thm} Let $\sigma$ be a proper primitive substitution. Then the $K$-groups of $B_{\sigma }$ are given by
\begin{equation*}
K_{0}(B_{\sigma })\simeq \mathbb{Z}\oplus \text{Coker}(1-\tau _{\sigma }), ~~~
K_{1}(B_{\sigma })\simeq \mathbb{Z}\oplus \text{Ker}(1-\tau _{\sigma }).
\end{equation*} 
\end{thm}
\begin{proof}
It follows from \cite{Cu2} and the Pimsner-Voiculescu sequence in \cite{Pi} 
that there is a exact sequence
$$
\begin{CD}
    {K_{\sigma }}@>{i_*-\rho _{\ast }}>> K_{\sigma }
@>i_*>> K_{0}(B_{\sigma }) \\
   @A{}AA
    @.
     @VV{}V \\
   K_{1}(B_{\sigma }) @<<i_*< \mathbb{Z} @<<{i_*-\rho _{\ast }}< \mathbb{Z}
\end{CD}
$$
where $\rho _{\ast }$ is the induced map of $\rho$ on $K_{0}(\mathcal{F}_{\sigma })$ and $K_{1}(\mathcal{F}_{\sigma })$. 
Now $\rho _{\ast }=\tau _{\sigma }$ on $K_{0}(\mathcal{F}_{\sigma })$ and 
$\rho _{\ast }=1$ on $K_{1}(\mathcal{F}_{\sigma })$ and we conclude the proof.
\end{proof}

\begin{example}
Let $\sigma$ be a para-periodic proper substitution. Since the map $\tau _{\sigma }$ 
on $K_{\sigma }$ is given by $\tau _{\sigma } (x+H_{\sigma })=\frac{x}{N_{\sigma }}+H_{\sigma }$, $K_{0}(B_{\sigma })\simeq \mathbb{Z}$ and $K_{1}(B_{\sigma })\simeq \mathbb{Z}$.
\end{example}

\end{document}